\documentclass[reqno,centertags,12pt]{amsart}

\usepackage{amssymb,amsmath,amsthm,dsfont}

\makeatletter \def\@strippedMR{} \def\@scanforMR#1#2#3\endscan{%
  \ifx#1M\ifx#2R\def\@strippedMR{#3}%
  \else\def\@strippedMR{#1#2#3}%
  \fi\fi} \renewcommand\MR[1]{\relax \ifhmode\unskip\spacefactor3000
  \space\fi \@scanforMR#1\endscan
  MR\MRhref{\@strippedMR}{\@strippedMR}} \makeatother

\newcommand{\R}{\mathbb{R}} \newcommand{\Z}{\mathbb{Z}}
\newcommand{\T}{\mathbb{T}} \newcommand{\C}{\mathbb{C}}
\newcommand{\N}{\mathbb{N}}\renewcommand{\S}{\mathbb{S}}

\theoremstyle{plain} \newtheorem{theorem}{Theorem}[section]
\newtheorem{lemma}[theorem]{Lemma}
\newtheorem{coro}[theorem]{Corollary}
\newtheorem{prop}[theorem]{Proposition}

\theoremstyle{definition} \newtheorem{definition}[theorem]{Definition}
\theoremstyle{remark} \newtheorem{remark}{Remark}

\DeclareMathOperator{\Id}{Id}

  \newcommand{\eps}{\varepsilon} \newcommand{\lb}{\langle}
\newcommand{\rb}{\rangle}
\newcommand{\ls}{\lesssim}\newcommand{\gs}{\gtrsim}

\begin{document}

\title[The quintic NLS on $3d$ Zoll manifolds]{The quintic nonlinear
  Schr\"odinger equation on three-dimensional Zoll manifolds}
\author[S.~Herr]{Sebastian~Herr}


\address{Universit\"at Bonn, Mathematisches Institut, Endenicher Allee
  60, 53115 Bonn, Germany}
\curraddr{Universit\"at Bielefeld, Fakult\"at f\"ur Mathematik,
  Postfach 100131, 33501 Bielefeld, Germany}
\email{herr@math.uni-bielefeld.de}

\begin{abstract}
  Let $(M,g)$ be a three-dimensional smooth compact Riemannian
  manifold such that all geodesics are simple and closed with a common
  minimal period, such as the $3$-sphere $\S^3$ with canonical metric.
  In this work the global well-posedness problem for the quintic
  nonlinear Schr\"odinger equation $i\partial_t u+\Delta u=\pm
  |u|^4u$, $u|_{t=0}=u_0$ is solved for small initial data $u_0$ in
  the energy space $H^1(M)$, which is the scaling-critical space.
  Further, local well-posedness for large data, as well as
  persistence of higher initial Sobolev regularity is obtained. This extends
  previous results of Burq-G\'erard-Tzvetkov to the endpoint case.
\end{abstract}
\maketitle
\section{Introduction and main result}\label{sect:intro_main}
\noindent
Let $(M,g)$ be a smooth Riemannian manifold without boundary and let
$\Delta=\Delta_g$ denote the (negative) Laplace-Beltrami operator on
$M$. Consider the Cauchy problem
\begin{equation}\label{eq:nls}
  \begin{split}
    i\partial_t u+\Delta u=&\pm |u|^4u\\
    u|_{t=0}=&\phi \in H^s(M)
  \end{split}
\end{equation}

If $u:(-T,T)\times M \to \C$ is a sufficiently nice solution to
\eqref{eq:nls} one easily verifies conservation of mass and energy
\begin{align}
  \label{eq:l2}
  \mathbf{m}(u(t))=&\tfrac12 \int_{M} |u(t,x)|^2 dx=\mathbf{m}(\phi),\\
  \label{eq:e}
  \mathbf{e}(u(t))=&\tfrac12 \int_{M} |\nabla u(t,x)|^2\pm \tfrac13
  |u(t,x)|^6 dx=\mathbf{e}(\phi).
\end{align}
Therefore, the Sobolev space $H^1(M)$ is the natural energy space for
\eqref{eq:nls}, in which for small initial data the local and the
global well-posedness problem are at the same level of difficulty.
Also, in the three-dimensional Euclidean case $(\R^3,\delta_{ij})$ the
scaling
\[u(t,x)\rightarrow \lambda^{\frac12} u(\lambda^2 t, \lambda x) \qquad
(\lambda>0)\] maps solutions onto solutions and does not alter the
$\dot{H}^1(\R^3)$-norm. Therefore, the energy space is called
(scaling-)critical.

In continuation of the line of research initiated by
Burq-G\'erard-Tzvetkov in \cite{BGT04,BGT05a,BGT05b,BGT07} we focus on
three-dimensional Zoll manifolds such as the sphere $M=\S^3$ with
canonical metric $g$, see \cite{Be78} for more information on the
geometric assumption, and \eqref{eq:spec} below. The sub-quintic
problem on $\S^3$ is solved in \cite[Theorem 1]{BGT05a}, and it is proved
that the super-quintic problem is ill-posed in \cite[Appendix
A]{BGT05a}. Well-posedness in $H^1(M)$ for the quintic nonlinear
Schr\"odinger equation in $\S^3$ is formulated as an open problem in
\cite[p. 257, l. 11]{BGT05a}, and it is shown in \cite{BGT07} that the
second Picard-iteration is bounded. This is the starting point for the
present paper, in which we prove

\begin{theorem}\label{thm:main} Let $s\geq 1$ and $(M,g)$ be a
  three-dimensional smooth compact Riemannian manifold such that all
  geodesics are simple and closed with a common minimal period. Then,
  the initial value problem \eqref{eq:nls} is locally well-posed in
  $H^s(M)$, and globally well-posed in $H^s(M)$ if the data is small
  in $H^1(M)$.
\end{theorem}
We refer the reader to Theorem \ref{thm:main-tech} in Section
\ref{sect:pf_main} for a more precise statement of the main result.
Theorem \ref{thm:main} completes the small data well-posedness theory
for $3d$ Zoll manifolds as we push it to the critical space $H^1(M)$, and Thomann's work \cite{Th08} shows that the problem is ill-posed in $H^s(M)$ for $s<1$ and analytic $(M,g)$.

Zoll manifolds have the property that the spectrum of the
Laplace-Beltrami operator $\Delta$ is clustered around a sequence of
squares \cite{DG,CdV,W77}, see \eqref{eq:spec} below. On the other
hand, the spectral cluster estimates \cite{BGT05a} are optimal on
spheres. This constitutes a sharp contrast to the case of the flat
rational torus $M=\T^3$ where we have recently established the
analogous result to Theorem \ref{thm:main}, see \cite{HTT10a}.

The choice of Zoll manifolds as our setup is motivated from
\cite{BGT04,BGT05a,BGT05b}. Also, one of our main ingredients in the
proof -- the trilinear spectral cluster estimates (Lemma
\ref{lem:tri-sogge}) -- are provided in \cite{BGT05a}. In this paper
we use critical function space techniques which have been introduced
by Tataru and Koch-Tataru \cite{KT05}, see also
\cite{HHK09} for further details. They have already been applied to related problems, namely
energy-critical Schr\"odinger equations on $M=\T^3$ \cite{HTT10a} and
$M=\R^2\times \T^2, \R^3\times \T$ \cite{HTT10b}.
On a technical level, however, we face different challenges in
this paper as the estimates for space and time variables will be
decoupled, and Galilean invariance and fine orthogonality arguments in
the spatial frequencies are unavailable. Actually, our function spaces
are slightly different and the general strategy of proof of Corollary
\ref{cor:tri-str-v2} also provides an alternate approach to
\cite{HTT10a}.

The paper is organized as follows: In Section \ref{sect:not_fs} we
describe the geometric and functional setup. Section \ref{sect:est}
starts by collecting known estimates on exponential sums and spectral
projectors, and after some preparation it concludes with the key estimate
of this work in Corollary \ref{cor:tri-str-v2}. Section
\ref{sect:pf_main} contains the main nonlinear estimate and a precise
statement of the main result in Theorem \ref{thm:main-tech}. In the
Appendix we describe the necessary modifications with respect to
Bourgain's paper \cite{B89} in order to conclude Lemma \ref{lem:lp}.

\section{Notation and function spaces}\label{sect:not_fs}
\noindent
Since $M$ is compact, the spectrum $\sigma(-\Delta)$ of
$-\Delta=-\Delta_g$ is discrete and we list the nonnegative
eigenvalues $0=\lambda_0^2\leq \lambda_1^2\leq \ldots \leq \lambda_n^2
\to +\infty$. Define
$h_k:L^2(M)\to L^2(M)$ to be the spectral projector onto the
eigenspace $E_k$ corresponding to the eigenvalue $\lambda_k^2$. We
have the orthogonal decomposition
\[
L^2(M)=\bigoplus\limits_{k=0}^\infty E_k,
\]
Following \cite{BGT04,BGT05a,BGT05b} we assume that $(M,g)$ is a
three-dimensional Zoll manifold, i.e. all geodesics are simple and
closed with a common minimal period $T$, and without loss we may
assume that $T=2\pi$.  In fact, we are using this assumption only to
conclude that the spectrum is clustered around the sequence $\mu_n^2$
where $\mu_n:=(n+\alpha/4)$ (for convenience of notation we define
$\mu_0:=0$). More precisely, there exist $\alpha,E\in \N$ such that
\begin{equation}\label{eq:spec} \sigma(-\Delta)\subset
  \cup_{n=1}^\infty I_n, \text{ where } I_n= [\mu_n^2-E,\mu_n^2+E],
\end{equation}
see \cite{DG,CdV,W77}.  By adding a bounded interval $I_0$, and
increasing $\alpha$ and relabeling if necessary, we may assume without
loss of generality that $\sigma(-\Delta)\subset \cup_{n=0}^\infty
I_n$, where
\begin{equation}\label{eq:spec-sp}
  I_n= [\mu_n^2-E,\mu_n^2+E] \; \text{ for } n \geq 1, \text{ and }I_0=[-B,B],
\end{equation}
with the additional property that all these intervals are pairwise
disjoint.  Let $p_n=\sum_{k \in \N_0: \lambda_k^2\in I_n} h_k$ for
$n\in \N_0$. Note that this is a spectral projector to which the
result of \cite[Theorem 3]{BGT05a} applies.  We have
\[\sum_{n=0}^\infty p_n=\Id.\] For subsets $J\subset \R$ we define
\[P_J=\sum_{n\in J\cap \N_0} p_n.\] Specifically, for dyadic numbers
$N=1,2,4,\ldots$ we write
\[
P_N=P_{[N,2N)}, \quad P_0=p_0, \text{ such that } \sum_{N \geq 0}
P_N=\Id,
\]
where we add up all dyadic $N\geq 1$ and $N=0$.

We define $H^s(M)=(1-\Delta_g)^{-s/2}L^2(M)$ with
\[
\|f\|_{H^s(M)}^2=\sum_{j=0}^\infty\lb \lambda_j \rb^{2s}
\|h_jf\|^2_{L^2(M)},
\]
where $\lb x\rb=(1+|x|^2)^{1/2}$ and observe that
\begin{equation*}\|f\|_{H^s(M)}^2 \sim \sum_{n=0}^\infty \lb n\rb
  ^{2s} \|p_{n} f\|^2_{L^2(M)}\sim \sum_{N \geq 0}\lb N \rb ^{2s}\|P_N
  f\|^2_{L^2(M)}.
\end{equation*}
\begin{remark}
  In the case $M=\S^3$ the eigenfunctions in $E_k=h_kL^2(\S^3)$ are
  precisely the spherical harmonics of degree $k$, with eigenvalues
  \[\lambda_k^2=(k+1)^2-1,\]
  and $\dim(E_k)=(k+1)^2$, see e.g.  \cite[Section 8.4]{TII}. With
  $E=1$, $\alpha=4$ in \eqref{eq:spec} we have $\lambda_k^2\in I_n$ if
  and only if $k=n$, and $\lambda_n^2=\mu_n^2-1$ for $n\geq 2$.
\end{remark}

For technical purposes we introduce the operator $\widetilde{\Delta}$
defined by
\[-\widetilde{\Delta}\phi = \sum_{n=1}^\infty \mu^2_n p_n \phi,\]
which is similar to the construction in \cite[formula (3.5)]{BGT05b}.

Let us quickly review the theory of the critical function spaces $U^p$
and $V^p$ which have been introduced in the context of dispersive PDEs
by Tataru and Koch-Tataru, see \cite{KT05}. We refer the reader to the
papers \cite{HHK09,HTT10a} for more details.
Let $\chi_I:\R\to\R$ denote the sharp characteristic function of a set
$I\subset \R$.  Let $\mathcal{Z}$ be the set of finite partitions
$-\infty<t_0<t_1<\ldots<t_K\leq \infty$ of the real line. If
$t_K=\infty$, we use the convention that $v(t_K):=0$ for all functions
$v: \R \to L^2$.

\begin{definition}\label{def:uv}
  Let $1\leq p <\infty$.
  \begin{enumerate}
  \item Any step-function $a:\R \to L^2$,
    \begin{equation*}
      a=\sum_{k=1}^K\chi_{[t_{k-1},t_k)}\phi_{k-1}
    \end{equation*}
with $\{t_k\} 
      \in \mathcal{Z}$,  $\{\phi_k\} \subset L^2$ s.t. $\sum_{k=0}^{K-1}\|\phi_k\|_{L^2}^p=1$,
is called a $U^p$-atom. We define $U^p$ as the corresponding atomic
space, i.e. the space of all $u: \R \to L^2$ which can be written as
\begin{equation}\label{eq:atomic-dec}
u=\sum_{j=1}^\infty \alpha_j a_j
\end{equation}
with $U^p$-atoms $a_j$, and $\{\alpha_j\}\in \ell^1(\N, \C)$. The
norm of a function $u\in U^p$ is defined as $\inf\sum_{j=1}^\infty |\alpha_j|$, where
the infimum is taken over all atomic representations
\eqref{eq:atomic-dec} of $u$.
  \item Define $V^p$ as the space of all right-continuous
    functions $v:\R\to L^2$ s.t.
    \begin{equation}\label{eq:norm_v}
      \|v\|_{V^p}:=\sup_{\{t_k\} \in \mathcal{Z}} \left(\sum_{k=1}^{K}
        \|v(t_{k})-v(t_{k-1})\|_{L^2}^p\right)^{\frac{1}{p}},
    \end{equation}
    is finite (here we use the convention
      $v(\infty)=0$), and additionally satisfying $\lim_{t\to
      -\infty}v(t)=0$.
  \item If $L^2=L^2(M;\C)$, and $A$ denotes either the standard
    Laplacian $\Delta$ or $\widetilde{\Delta}$, we also define
    $U^p_A=e^{itA}U^p$ and $V^p_A=e^{itA}V^p$.
  \end{enumerate}
\end{definition}

\begin{remark}\label{rmk:prop_uv}
  \begin{enumerate}
\item\label{it:hhk} Note that the space $V^p$ corresponds to $V^p_{-,rc}$ in \cite{HHK09}.
  \item\label{it:prop_uv1} The spaces $U^p,V^p$ and $U^p_A,V^p_A$ are
    Banach spaces of bounded functions which are right-continuous and
    tend to $0$ as $t\to -\infty$.
  \item\label{it:prop_uv2} For $1\leq p<q <\infty$ it holds
    \[U^p_A\hookrightarrow V^p_A\hookrightarrow U^q_A \hookrightarrow
    L^\infty(\R;L^2).\]
  \end{enumerate}
\end{remark}

Our aim is to control the evolution up to time $T\sim 1$. But on
bounded time intervals the flows associated to the operators $\Delta$
and $\widetilde{\Delta}$ stay close, a statement which is made precise
next.
\begin{lemma}\label{lem:mod_sp}
  Let $1\leq p<\infty$ and $\tau$ be a bounded time interval, and let
  $u:\R\to L^2(M;\C)$ be supported in $\tau$.  Then, $u \in
  U^p_\Delta$ if and only if $u \in U^p_{\widetilde{\Delta}} $, with
  equivalent norms.
\end{lemma}
\begin{proof}
  Let $\psi\in C_0^\infty (\R)$ be a smooth cutoff function which is
  constantly equal to $1$ on $\tau$. The claim follows from the fact
  that
  \[\psi e^{\pm it(\Delta-\widetilde{\Delta})}:U^p\to U^p\]
  are bounded operators. It suffices to consider an atom $a$. With
  $B_j:=\sum_{n=0}^\infty \sum_{k \in \N_0: \lambda_k^2\in I_n } (\pm
  \lambda_k^2\mp \mu_n^2)^j h_k$ we write
  \[\psi(t) e^{\pm it(\Delta-\widetilde{\Delta})}
  a(t)=\sum_{j=0}^\infty \frac{\psi(t) (it)^j}{j!} B_j a(t),\] In view
  of \eqref{eq:spec-sp} we have $\|B_j\|_{L^2\to L^2}\leq b^j$, where
  $b=\max\{E,B\}\geq 1$.  This implies $B_j a\in U^p$ with $\|B_j
  a\|_{U^p}\leq b^j$.  Also, multiplication by $\psi$, as well as
  multiplication by $t$ on the support of $\psi$ are bounded
  operations in $U^p$, which follows by duality \cite[Remark
  2.11]{HHK09}. This implies
  \[
  \|t^j \psi B_ja\|_{U^p}\leq c^j.
  \]
  The claim follows, since $\sum_{j=0}^\infty \frac{c^j}{j!}<+\infty$
  and $U^p$ is a complete space.
\end{proof}

\begin{definition}\label{def:xy} Let $s\in \R$.
  \begin{enumerate}
  \item We define $X^s$ as the space of all functions $u: \R \to
    H^s(M; \C)$ such that the maps $P_N(u(\cdot)):\R\to L^2(M;\C)$ are
    in $U^2_\Delta$ for all dyadic $N\geq 0$, endowed with the norm
    \begin{equation}\label{eq:xsnorm}
      \|u\|_{X^s}:=\left(\sum_{N\geq 0} \lb N \rb^{2s} \|P_N u\|_{U^2_\Delta}^2\right)^{\frac12}.
    \end{equation}
  \item We define $Y^s$ as the space of all functions $u:\R \to
    H^s(M;\C)$ such that the maps $P_N(u(\cdot)):\R\to L^2(M;\C)$ are
    in $V^2_\Delta$ for all dyadic $N\geq 0$, equipped with the norm
    \begin{equation}\label{eq:ysnorm}
      \|u\|_{Y^s}
      :=\left(\sum_{N \geq0}
        \lb N \rb^{2s} \|P_N u\|_{V^2_\Delta}^2\right)^{\frac12}.
    \end{equation}
  \end{enumerate}
  For a time interval $\tau \subset \R$ we define $X^s(\tau)$ and
  $Y^s(\tau)$ to be the corresponding restriction space.
\end{definition}

The above remark implies
\[
U^2_\Delta \hookrightarrow X^0 \hookrightarrow Y^0 \hookrightarrow
V^2_\Delta.
\]
Moreover, there is a useful interpolation type property of $U^p$ and
$V^p$ spaces, see \cite[Proposition 2.20]{HHK09} and \cite[Lemma
2.4]{HTT10a}.
\begin{lemma}\label{lem:interpol}
  Let $q_1,q_2,q_3>2$, $\tau$ be a time interval, and
  \[T:U_\Delta^{q_1}\times U^{q_2}_\Delta\times U^{q_3}_\Delta\to
  L^2(\tau \times M)\] be a bounded, tri-linear operator with
  $\|T(u_1,u_2,u_3)\|_{L^2} \leq C
  \prod_{j=1}^3\|u_j\|_{U^{q_j}_\Delta}$.  In addition, assume that
  there exists $C_2\in (0,C]$ such that the estimate
  $\|T(u_1,u_2,u_3)\|_{L^2} \leq C_2
  \prod_{j=1}^3\|u_j\|_{U^2_\Delta}$ holds true.  Then, $T$ satisfies
  the estimate
  \[
  \|T(u_1,u_2,u_3)\|_{L^2} \ls {} C_2
  (\ln\frac{C}{C_2}+1)^3\prod_{j=1}^3\|u_j\|_{V^2_\Delta}, \quad u_j
  \in V^2_\Delta,\; j=1,2,3.
  \]
\end{lemma}

Let $\tau=[a,b)$, $f\in L^1(\tau;L^2(M;\C))$ and define
\begin{equation}\label{eq:duhamel}
  \mathcal{I}(f)(t):=\int_{a}^t e^{i(t-s)\Delta} f(s) ds
\end{equation}
for $t \in \tau$ and $\mathcal{I}(f)(t)=0$ for $t<a$ and
$\mathcal{I}(f)(t)=\mathcal{I}(f)(b)$ for $t\geq b$.

\begin{lemma}\label{lem:linear}
  Let $s\in \R$, $\tau=[a,b)\subset \R$.
  \begin{enumerate}
  \item For all $u_0 \in H^s(M)$ and
    $u(t):=\chi_{\tau}(t)e^{it\Delta}u_0$ we have $u \in X^s(\tau)$
    and
    \begin{equation}\label{eq:linear1}
      \|u\|_{X^s(\tau)}\ls {} \|u_0\|_{H^s}.
    \end{equation}
  \item Let $P_N f \in L^1(\tau;L^2(M))$ for all $N\geq0$. Then,
    $\sum_{N\geq0}\mathcal{I}(P_Nf)=:\mathcal{I}(f)$ converges in
    $X^s(\tau)$ and
    \begin{equation}\label{eq:linear2}
      \|\mathcal{I}(f)\|_{X^s(\tau)}
      \leq \sup_{\|v\|_{Y^{-s}(\tau)}=1} \Big|\sum_{N\geq 0}\int_\tau \int_{M}
      P_Nf(t,x) \overline{v(t,x)}
      dx dt \Big|,
    \end{equation}
    provided that the r.h.s in \eqref{eq:linear2} is finite.
  \end{enumerate}
\end{lemma}
We refer the reader to \cite[Propositions 2.10 and 2.11]{HTT10a} and
\cite[Propositions 2.8 and 2.10]{HHK09} for analogous
statements and proofs, which apply here with trivial modifications.

\section{Linear and multilinear estimates}\label{sect:est}
We use an extension of Bourgain's estimate \cite[Proposition 1.10 and
Section 4]{B89} on exponential sums which is due to
Burq--G\'erard--Tzvetkov \cite[Lemma 5.3]{BGT07} in the case $p=6$,
$\mu_n=n$. For convenience we choose $\tau_0=[0,32\pi]$ as our base
time interval, as this is a joint period of $e^{-it \mu_n^2}$.
\begin{lemma}
  \label{lem:lp}
  Let $p>4$ and $\alpha\in \N_0$. It holds that
  \begin{equation}\label{eq:lp}
    \Big\|\sum_{n \in \Z \cap J} c_n e^{-it \mu_n^2}\Big\|_{L^p_t(\tau_0)}\ls {}
    N^{1/2-2/p}\Big(\sum_{n \in \Z\cap J} |c_n|^2\Big)^{\frac12}
  \end{equation}
  for every $J=[b,b+N]$ with $N\geq 1$ and the sequence
  $\mu_n=n+\alpha/4$.
\end{lemma}
A proof can be found in Appendix \ref{asect:p}.  We will also rely on
the trilinear spectral projector estimate of Burq--G\'erard--Tzvetkov
\cite{BGT05a}, which is valid on every smooth Riemannian three-manifold
$(M,g)$.
\begin{lemma}[\cite{BGT05a}, Theorem 3]
  \label{lem:tri-sogge} Let $0<\epsilon\ll 1$. For all integers
  $n_1\geq n_2 \geq n_3 \geq 0$ and $f_1,f_2,f_3\in L^2(M)$ the
  estimate
  \begin{equation}\label{eq:tri-sogge}
    \|p_{n_1}f_1 p_{n_2}f_2 p_{n_3}f_3\|_{L^2(M)}
    \ls {} \lb n_2\rb ^{1/2+\epsilon}\lb n_3\rb ^{1-\epsilon}
    \prod_{k=1}^3\|p_{n_k } f_k\|_{L^2(M)}
  \end{equation}
  holds true.
\end{lemma}

In the nonlinear analysis we need another useful and well-known
estimate concerning the spectral localization of products of
eigenfunctions.
\begin{lemma}\label{lem:decay}
  If $N_0\gg N_1,N_2,N_3$ are dyadic and $\gamma \geq 1$, then
  \begin{equation}\label{eq:decay}
    \Big|\int_M P_{N_0} f_0P_{N_1} f_1P_{N_2} f_2P_{N_3} f_3 dx\Big|\ls {} N_0^{-\gamma} \prod_{j=0}^3 \|P_{N_j} f_j\|_{L^2(M)}
  \end{equation}
  where the implicit constant depends only on $\gamma$.
\end{lemma}
\begin{proof}
  For single eigenfunctions it can be found in \cite[Lemma
  2.6]{BGT05b} (written for $d=2$) or more generally in \cite[Section
  4]{H10}. By the Weyl asymptotic the number of eigenvalues
  $\lambda_k^2 \in I_{n_j}$ grows at most like $n_j^2$ in dimension
  $d=3$, which implies
  \begin{align*}
    & \Big|\int_M P_{N_0} f_0P_{N_1} f_1P_{N_2} f_2P_{N_3} f_3 dx\Big|\\
    \leq {}& \sum_{n_j\sim N_j}\sum_{k_j: \lambda_{k_j}^2\in I_{n_j}^2}\Big|\int_M h_{k_0} f_0h_{k_1} f_1h_{k_2} f_2h_{k_3} f_3 dx\Big|\\
    \ls {} & \sum_{n_j\sim N_j}\sum_{k_j: \lambda_{k_j}^2\in
      I_{n_j}^2} k_0^{-\gamma-10} \prod_{j=0}^3 \|h_{k_j} f_j\|_{L^2}
    \ls {} N_0^{-\gamma} \prod_{j=0}^3 \|P_{N_j} f_j\|_{L^2},
  \end{align*}
  by Cauchy-Schwarz, cp. also \cite[Lemma 2.7]{BGT05b} for similar
  arguments.
\end{proof}
Note that Lemma \ref{lem:decay} is trivial in specific cases such as
$M=\S^3$ with canonical metric.

For later reference we explicitely state a crude bound which
disregards all oscillations in time.
\begin{lemma}\label{lem:crude}
  For all $u_1,u_2,u_3\in L^\infty(\tau;L^2(M))$, and dyadic $N_1\geq
  N_2\geq N_3\geq 0$ and time intervals $\tau$ the estimate
  \begin{equation}\label{eq:crude}
    \|P_{N_1}u_1P_{N_2}u_2P_{N_3}u_3\|_{L^2(\tau\times M)}\\
    \ls {}  |\tau|^{\frac12} N_2^{\frac32} N_3^{\frac32}\prod_{j=1}^3\|u_j\|_{L^\infty(\tau;L^2(M))}
  \end{equation}
  holds true.
\end{lemma}
\begin{proof}
  By H\"older's inequality we obtain
  \begin{align*}
    &\|P_{N_1}u_1P_{N_2}u_2P_{N_3}u_3\|_{L^2(\tau\times M)}\\
    \leq & |\tau|^{\frac12} \|P_{N_1}u_1\|_{L^\infty(\tau;L^2(M)}
    \|P_{N_2}u_2\|_{L^\infty(\tau \times
      M)}\|P_{N_3}u_3\|_{L^\infty(\tau \times M)}.
  \end{align*}
  For $t$ fixed we have
  \begin{align*}
    &\|P_{N_j}u_j(t)\|_{L^\infty(M)}\leq \sum_{n_j \sim N_j}\|p_{n_j}u_j(t)\|_{L^\infty(M)}\\
    \ls {}& \sum_{n_j \sim N_j}n_j \|p_{n_j}u_j(t)\|_{L^2(M)}\ls {}
    N_j^{\frac32} \|P_{N_j}u_j(t)\|_{L^2(M)}
  \end{align*}
  by Sogge's estimate \cite[Proposition 2.1]{S88} and the
  Cauchy-Schwarz inequality.
  The claim follows by taking the supremum in $t$.
\end{proof}

Next, we prove (dyadic) Strichartz estimates in a restricted range,
generalizing \cite[Proposition 5.1]{BGT07}.
\begin{lemma}\label{lem:str}
  Let $p>4$. Then, for all $N\geq 0$ we have
  \begin{equation}\label{eq:str}
    \|P_N e^{it\Delta}\phi\|_{L^p(\tau_0\times M)}\ls {}
    \lb N\rb^{\frac32-\frac5p}\|\phi\|_{L^2}.
  \end{equation}
  and
  \begin{equation}\label{eq:str-up}
    \|P_N u\|_{L^p(\tau_0\times M)}\ls {} \lb N\rb^{\frac32-\frac5p}\|u\|_{U^p_\Delta}.
  \end{equation}
\end{lemma}
\begin{proof} {\it a)} First, we prove estimate \eqref{eq:str} with
  $\widetilde{\Delta}$ replacing $\Delta$. We write
  \[P_N e^{it\widetilde{\Delta}}\phi(x)=\sum_{n\sim
    N}e^{-it\mu_n^2}p_n\phi(x),\] and Lemma \ref{lem:lp} yields
  \[
  \|P_N e^{it\widetilde{\Delta}}\phi(x)\|_{L^p(\tau_0)}\ls {} \lb
  N\rb^{\frac12-\frac2p}\Big(\sum_{n\sim
    N}|p_n\phi(x)|^2\Big)^{\frac12}
  \]
  Integration in $x$, an application of Minkowski's inequality and the
  dual of Sogge's estimate \cite[formula (2.3)]{S88} imply that
  \begin{align*}
    &\Big\|\Big(\sum_{n\sim N}|p_n\phi(x)|^2\Big)^{\frac12}\Big\|_{L^p(M)}\ls {} \Big(\sum_{n\sim N}\|p_n\phi\|_{L^p}^2\Big)^{\frac12}\\
    \ls {} & \Big(\sum_{n\sim N} \lb
    n\rb^{2-\frac6p}\|p_n\phi\|_{L^2}^2\Big)^{\frac12}\ls {} \lb
    N\rb^{1-\frac3p}\|\phi\|_{L^2},
  \end{align*}
  where we have used orthogonality in the last step.

  {\it b)} Now, let $u\in U^p_\Delta$. By Lemma \ref{lem:mod_sp} it
  suffices to prove the bound \eqref{eq:str-up} for a
  $U^p_{\widetilde{\Delta}}$-atom
  \[u(t,x)=\sum_{k=1}^{K}\chi_{[t_{k-1},t_{k})}e^{it\widetilde{\Delta}}\phi_{k-1}(x),
  \quad \sum_{k=0}^{K-1}\|\phi_{k}\|^p_{L^2}=1.\] Estimate
  \eqref{eq:str} yields
  \begin{align*}
    \|P_N u\|_{L^p(\tau_0\times M)}\leq {} & \Big(\sum_{k=1}^{K} \|P_N
    e^{it\widetilde{\Delta}}\phi_{k-1}\|^p_{L^p(\tau_0\times
      M)}\Big)^{\frac1p}\\
    \ls {} & \lb N\rb^{\frac32-\frac5p} \Big(\sum_{k=1}^{K}
    \|P_N\phi_{k-1}\|^p_{L^2(M)}\Big)^{\frac1p}\\
    \ls {} & \lb N\rb^{\frac32-\frac5p}.
  \end{align*}
  This proves the second bound \eqref{eq:str-up}, which in turn
  implies \eqref{eq:str} and the proof is complete.
\end{proof}
E.g. on the sphere $M=\S^3$, the restriction to dyadic frequency bands
can be removed by Littlewood-Paley theory \cite{St72}, but we do not
need it here. Also, the loss of derivatives precisely matches the loss
on $\R^3$ coming from the sharp Strichartz estimate and the Sobolev
embedding, and also Bourgain's bound on $\T^3$ \cite{B93a}.

For an interval $J$ we write $P_{N,J}=P_J P_N$. The next Proposition is an extension and
improvement of \cite[Theorem 5.1]{BGT07}.
\begin{prop}\label{prop:u2-tri-str}
  Let $\delta \in (0,\tfrac12)$ and $\eta>0$. Then, for all
  $u_1,u_2,u_3\in U^2_\Delta$, and dyadic $N_1\geq N_2\geq N_3\geq 0$
  the estimate
  \begin{equation}\label{eq:tri-str-u2}
    \begin{split}
      &\|P_{N_1}u_1P_{N_2}u_2P_{N_3}u_3\|_{L^2(\tau_0 \times M)}\\
      \ls {} & \Big(\frac{\lb N_2\rb}{\lb N_1\rb }+\frac{1}{\lb N_2
        \rb}\Big)^{\delta} \lb N_{2}\rb^{\frac12+\eta+\delta}\lb
      N_{3}\rb^{\frac32-\eta-\delta}\prod_{j=1}^3\|u_j\|_{U^2_\Delta}.
    \end{split}
  \end{equation}
  holds true.
\end{prop}
\begin{proof} {\it a)} In the case $N_3\leq N_2\leq 1$ the l.h.s. is
  bounded by our crude estimate \eqref{eq:crude} in conjunction with
  Remark \ref{rmk:prop_uv} \ref{it:prop_uv2}. Henceforth we assume
  $N_2\geq 1$. Since $\tau_0$ is a compact interval, it suffices to
  prove the corresponding bound in $U^2_{\widetilde{\Delta}}\times
  U^2_{\widetilde{\Delta}}\times U^2_{\widetilde{\Delta}}$ by Lemma
  \ref{lem:mod_sp}. Further, by definition of the spaces it suffices
  to consider $U^2_{\widetilde{\Delta}}$-atoms
  \[
  u_j(t,x)=\sum_{k_j=1}^{K_j}\chi_{[t_{k_j-1}^{(j)},t_{k_j}^{(j)})}e^{it\Delta}\phi_{k_j-1}^{(j)}(x)
  , \quad \sum_{k_j=0}^{K_j-1}\|\phi_{k_j}^{(j)}\|^2_{L^2}=1,
  \]
  for $j=1,2,3$. However, in this case
  \[
  \|\prod_{j=1}^3P_{N_j}u_j\|^2_{L^2(\tau_0 \times M)}\leq
  \sum_{k_1,k_2,k_3}\|\prod_{j=1}^3
  P_{N_j}e^{it\Delta}\phi_{k_j-1}^{(j)}\|^2_{L^2(\tau_0 \times M)},\]
  so the claim \eqref{eq:tri-str-u2} follows if we can show it in the
  case
  \begin{equation}\label{eq:lin}
    u_j=e^{it\widetilde{\Delta}}\phi_j, \quad j=1,2,3.
  \end{equation}

  {\it b)} Assume \eqref{eq:lin}. We define the partition
  \[
  \N_0=\dot\cup_{m \in \N_0}J_m\quad \text{ where } J_m=[m
  N_2^2/N_1,(m+1)N_2^2/N_1)\cap \N_0
  \]
  in order to proceed similarly to \cite[proof of (26),
  p. 341--342]{HTT10a} and \cite{HTT10b}: For fixed $x\in M$
  it holds
  \[
  \|P_{N_1}u_1P_{N_2}u_2P_{N_3}u_3(x)\|_{L^2(\tau_0)}^2 \sim \sum_m
  \|P_{N_1,J_m}u_1P_{N_2}u_2P_{N_3}u_3(x)\|_{L^2(\tau_0)}^2,
  \]
due to almost orthogonality induced by the time oscillations.
Indeed, it holds that
  \begin{align*}
    &\lb P_{N_1,J_m}u_1P_{N_2}u_2P_{N_3}u_3(x),P_{N_1,J_{m'}}u_1P_{N_2}u_2P_{N_3}u_3(x)\rb_{L^2(\tau_0)}\\
    =&\sum_{n_j,{n_j}'\sim N_j \atop n_1\in J_{m}, {n_1}' \in J_{m'}}
    I_{n_1,n_2,n_3}^{{n_1}',{n_2}',{n_3}'}\prod_{j=1}^3
    p_{n_j}\phi_j(x) \overline{p_{{n_j}'}\phi_j(x)}
  \end{align*}
  where
  \begin{equation*}
    I_{n_1,n_2,n_3}^{{n_1}',{n_2}',{n_3}'}=\int_{\tau_0} e^{-it\mu_{n_1}^2}e^{-it\mu_{n_2}^2}e^{-it\mu_{n_3}^2}e^{it\mu_{{n_1}'}^2}e^{it\mu_{{n_2}'}^2}e^{it\mu_{{n_3}'}^2}dt.
  \end{equation*}
  For every $|m-m'|\gg 1$ and $n_1,{n_1}'\sim N_1$ such that $n_1 \in
  J_m, {n_1}' \in J_{m'}$ and all $n_2,{n_2}'\sim N_2$,
  $n_3,{n_3}'\sim N_3$ we have the following estimate for the phase
  \[
  |\sum_{j=1}^3(\mu_{{n_j}'}^2-\mu_{n_j}^2)|\geq
  |\mu_{{n_1}'}^2-\mu_{n_1}^2|-8N_2^2\gs |m-m'|N_2^2,
  \]
  because $\mu_{{n_1}'}+\mu_{n_1}\geq N_1$, which implies
  \[I_{n_1,n_2,n_3}^{{n_1}',{n_2}',{n_3}'}=0.\]

  {\it c)} By parts a) and b) the claim is reduced to showing that
  \begin{equation}\label{eq:lin-tri-str}
    \|P_{N_1,J}u_1P_{N_2}u_2P_{N_3}u_3\|_{L^2(\tau_0 \times M)}
    \ls {}  |J|^{\delta}  N_{2}^{\frac12+\eta}\lb N_{3}\rb
    ^{\frac32-\eta-\delta}\prod_{j=1}^3\|\phi_j\|_{L^2}.
  \end{equation}
  for $u_j$ of the form \eqref{eq:lin}, and intervals of length
  $|J|\geq 1$.  This is a refinement of \cite[Theorem 5.1]{BGT07},
  which is proved as follows: For fixed $x\in M$ we obtain by
  H\"older's inequality
  \begin{align*}
    &\|P_{N_1,J}u_1P_{N_2}u_2P_{N_3}u_3(x)\|_{L^2(\tau_0)}\\
    \leq {}
    &\|P_{N_1,J}u_1(x)\|_{L^{p_1}(\tau_0)}\|P_{N_2}u_2(x)\|_{L^{p_2}(\tau_0)}\|P_{N_3}u_3(x)\|_{L^{p_3}(\tau_0)}
  \end{align*}
  where $1/p_1+1/p_2+1/p_3=1/2$ which we choose to satisfy
  $4<p_1,p_2,p_3<+\infty$.  An application of \eqref{eq:lp} gives
  \begin{align*}
    \|P_{N_1,J}u_1(x)\|_{L^{p_1}_t(\tau_0)} \ls {} &
    |J|^{\frac12-\frac2{p_1}}\Big(\sum_{n_1\sim N_1 \atop n_1 \in J}
    |p_{n_1}\phi_1(x)|^2\Big)^{\frac12}
  \end{align*}
  and also
  \begin{align*}
    \|P_{N_j}u_j(x)\|_{L^{p_j}_t(\tau_0)} \ls {} &
    N_j^{\frac12-\frac2{p_j}}\Big(\sum_{n_j\sim N_j}
    |p_{n_j}\phi_j(x)|^2\Big)^{\frac12}
  \end{align*}
  for $j=2,3$. By integration with respect to $x\in M$ we obtain
  \begin{align*}
    &\|P_{N_1,J}u_1 P_{N_2}u_2 P_{N_3}u_3\|_{L^2(\tau_0\times M)}\\
    \ls {}&
    |J|^{\frac12-\frac2{p_1}}N_2^{\frac12-\frac2{p_2}}N_3^{\frac12-\frac2{p_3}}\Big(
    \sum_{n_j\sim N_j}
    \|\prod_{j=1}^3 p_{n_j}\phi_{j}\|_{L^2(M)}^2\Big)^{\frac12}\\
    \ls {}&
    |J|^{\frac12-\frac2{p_1}}N_2^{1+\epsilon-\frac2{p_2}}N_3^{\frac32-\epsilon-\frac2{p_3}}\prod_{j=1}^3\|\phi_{j}\|_{L^2(M)}
  \end{align*}
  for any small $\epsilon>0$, where we have used the trilinear
  spectral cluster estimate \eqref{eq:tri-sogge} in the last step,
  similar to the proof of \cite[Theorem 5.1]{BGT07}. The
  claim follows with $\delta=\tfrac12-\tfrac2{p_1}\in (0,\tfrac12)$ by
  choosing $p_2>4$ and $\epsilon>0$ small enough such that
  $\epsilon+\tfrac12-\tfrac{2}{p_2}= \eta$.
\end{proof}
Finally, we transfer the bound to $V^2_\Delta$ by interpolation and
obtain a result which corresponds to \cite[Propositon 3.5]{HTT10a} in
the case of $M=\T^3$. The argument, however, is slightly different from the one in \cite{HTT10a,HTT10b} as it
does not involve finer than dyadic scales.
\begin{coro}\label{cor:tri-str-v2}
  There exists $\alpha>0$, such that for all $u_1,u_2,u_3\in
  V^2_\Delta$, and dyadic $N_1\geq N_2\geq N_3\geq 0$ the estimate
  \begin{equation}\label{eq:tri-str-v2}
    \begin{split}
      &\|P_{N_1}u_1P_{N_2}u_2P_{N_3}u_3\|_{L^2(\tau_0 \times M)}\\
      \ls {} & \max\Big\{\frac{\lb N_3\rb}{\lb N_1\rb },\frac{1}{\lb
        N_2 \rb}\Big\}^{\alpha} \lb N_{2}\rb \lb
      N_{3}\rb\prod_{j=1}^3\|u_j\|_{V^2_\Delta}.
    \end{split}
  \end{equation}
  holds true.
\end{coro}
\begin{proof}
  We restrict our attention to the nontrivial case $N_1\geq 1$, and
  treat the two cases
  \[\mathrm{a)}\, N_2^2\geq N_1 \qquad \mathrm{b)}\, N_2^2< N_1\]
  separately.

  {\it Case a)} For $p,q>4$ satisfying $\tfrac2p+\tfrac1q=\tfrac12$ we
  exploit \eqref{eq:str-up}
  \begin{align*}
    &\|P_{N_1}u_1P_{N_2}u_2P_{N_3}u_3\|_{L^2(\tau_0 \times M)}\\
    \leq {}&\|P_{N_1}u_1\|_{L^p(\tau_0 \times
      M)}\|P_{N_2}u_2\|_{L^p(\tau_0
      \times M)}\|P_{N_3}u_3\|_{L^q(\tau_0 \times M)}\\
    \ls {}&N_1^{\frac32-\frac5p} N_2^{\frac32-\frac5p}\lb
    N_3\rb^{\frac32-\frac5q}\|P_{N_1}u_1\|_{U^p_{\Delta}}\|P_{N_2}u_2\|_{U^p_{\Delta}}\|P_{N_3}u_3\|_{U^q_{\Delta}}.
  \end{align*}
  Let $\rho>0$ be small. We choose $p>4$ such that
  $\tfrac32-\tfrac{5}{p}=\tfrac14+\rho$. Then,
  \begin{equation}\label{eq:tri-str-up}
    \begin{split}
      &\|P_{N_1}u_1P_{N_2}u_2P_{N_3}u_3\|_{L^2(\tau_0 \times M)}\\
      \leq {}&\left(\frac{N_1}{N_2}\right)^{\frac14+\rho}
      N_2^{\frac12+2\rho}\lb
      N_3\rb^{\frac32-2\rho}\|P_{N_1}u_1\|_{U^p_{\Delta}}\|P_{N_2}u_2\|_{U^p_{\Delta}}\|P_{N_3}u_3\|_{U^q_{\Delta}}.
    \end{split}
  \end{equation}
  Interpolating \eqref{eq:tri-str-u2} and \eqref{eq:tri-str-up} via
  Lemma \ref{lem:interpol} yields
  \begin{align*}
    &\|P_{N_1}u_1P_{N_2}u_2P_{N_3}u_3\|_{L^2(\tau_0 \times M)}\\
    \ls {} & \Big(\frac{N_2}{ N_1 }\Big)^{\delta'} \lb
    N_{2}\rb^{\frac12+2\delta'}\lb
    N_{3}\rb^{\frac32-2\delta'}\prod_{j=1}^3\|u_j\|_{V^2_\Delta}.
  \end{align*}
  for small $\delta'>0$, because in the present case we have
  $N_2^2\geq N_1$.

  {\it Case b)} In this case where $N_2^2< N_1$, the key is to observe
  that \eqref{eq:tri-str-u2} provides the subcritical bound
  \begin{equation}\label{eq:tri-str-u2-sub}
    \|P_{N_1}u_1P_{N_2}u_2P_{N_3}u_3\|_{L^2(\tau_0 \times M)}
    \ls {} \lb N_{2}\rb^{\frac12+\eta}\lb
    N_{3}\rb\prod_{j=1}^3\|u_j\|_{U^2_\Delta}.
  \end{equation}
  for any $\eta>0$. On the other hand, estimate \eqref{eq:crude} and
  Remark \ref{rmk:prop_uv} \ref{it:prop_uv2} imply
  \begin{equation*}
    \|P_{N_1}u_1P_{N_2}u_2P_{N_3}u_3\|_{L^2(\tau_0 \times M)}
    \ls {} \lb N_{2}\rb^{\frac32}\lb
    N_{3}\rb^{\frac32}\prod_{j=1}^3\|u_j\|_{U^p_\Delta}.
  \end{equation*}
  for any $p\in [1,\infty)$, so interpolation via Lemma
  \ref{lem:interpol} yields
  \begin{equation}
    \|P_{N_1}u_1P_{N_2}u_2P_{N_3}u_3\|_{L^2(\tau_0 \times M)}
    \ls {} \lb N_{2}\rb^{\frac12+\eta'}\lb
    N_{3}\rb\prod_{j=1}^3\|u_j\|_{V^2_\Delta}.
  \end{equation}
  for any $\eta'>0$, which implies \eqref{eq:tri-str-v2} in this case.
\end{proof}

\section{The main result}\label{sect:pf_main}
\noindent
As usual, we rewrite the initial value problem as an integral equation
\begin{equation}\label{eq:int-nls}
  u(t)=e^{it\Delta}\phi\mp i \mathcal{I}(|u|^4u)(t).
\end{equation}

Now, we restate Theorem \ref{thm:main} in a more precise form.  Let us
denote the ball in $H^1(M)$ with center $\phi$ and radius $\eps$ by
$B_\eps(\phi)$.
\begin{theorem}\label{thm:main-tech}
  Let $(M,g)$ be a three-dimensional compact smooth Riemannian
  manifold with Laplace-Beltrami operator $\Delta$ satisfying the
  spectral condition \eqref{eq:spec}, and let $s\geq 1$.
  \begin{enumerate}
  \item\label{it:lwp} {\it (Local well-posedness)} For every
    $\phi_{\ast}\in H^1(M)$ there exists $\eps>0$ and
    $T=T(\phi_{\ast})>0$ such that the following holds true:
    \begin{enumerate}
    \item\label{it:a} For all initial data $\phi\in
      B_{\eps}(\phi_\ast)\cap H^s(M)$ the Cauchy problem
      \eqref{eq:int-nls} has a unique solution \[u=:\Phi(\phi) \in
      C([0,T);H^s(M)) \cap X^s([0,T)).\]
    \item\label{it:b} The solution constructed in Part \ref{it:a}
      obeys the conservation laws \eqref{eq:e} and \eqref{eq:l2}, and
      the flow map
      \[
      \Phi: B_{\eps}(\phi_\ast)\cap H^s(M)\to C([0,T);H^s(M)) \cap
      X^s([0,T))
      \]
      is Lipschitz continuous.
    \end{enumerate}
  \item\label{it:gwp} {\it (Global well-posedness for small data)}
    With $\phi_\ast=0$ there exists $\eps_0>0$ such that for all $T>0$
    the assertions \ref{it:a} and \ref{it:b} above hold true.
  \end{enumerate}
\end{theorem}

The proof is very similar to the proof of \cite[Theorems 1.1 and
1.2]{HTT10a}, as it is based on the following Proposition
\ref{prop:six-est} which corresponds to \cite[Proposition
4.1]{HTT10a}.

\begin{prop}\label{prop:six-est}
  Let $s\geq 1$. Then, for all intervals $\tau \subset \tau_0$ and all
  $u_j \in X^s(\tau)$, $j=1,\ldots, 5$, the estimate
  \begin{equation}\label{eq:six-est}
    \Big\|\mathcal{I}\Big(\prod_{j=1}^5\widetilde{u_j}\Big)\Big\|_{X^s(\tau)}\ls {}\sum_{k=1}^5
    \|u_k\|_{X^s(\tau)}\prod_{j=1; j\not=k}^5\|u_j\|_{X^1(\tau)}
  \end{equation}
  holds true, where $\widetilde{u_j}$ denotes either $u_j$ or its
  complex conjugate $\overline{u}_j$.
\end{prop}
\begin{proof} The proof is a variation of the proof of
  \cite[Proposition 4.1]{HTT10a}, so we focus on the new aspects here:
  Lemma \ref{lem:linear} implies that $\mathcal{I}(\prod_{j=1}^5
  \widetilde{u}_j)\in X^s(\tau)$ and
  \[
  \Big\|\mathcal{I}(\prod_{j=1}^5
  \widetilde{u}_j)\Big\|_{X^s(\tau)}\leq
  \sup_{\|u_0\|_{Y^{-s}(\tau)}=1} \Big|\sum_{N_0\geq 0}\int_\tau
  \int_{M}P_{N_0}\prod_{j=1}^5 \widetilde{u}_j \, \overline{u_0} dx
  dt\Big|,
  \]
  provided that the r.\ h.\ s.\ is finite.  Thus, by choosing suitable
  extensions (which we also denote by $u_j$) the claim is reduced to
  proving
  \begin{equation}\label{eq:6lin}
    \Big| \sum_{N_0\geq 0}\int_{\tau}\int_{ M} P_{N_0}\widetilde{u_0}\prod_{j=1}^5 \widetilde{u_j} \ dx dt\Big|
    \lesssim  \| u_0\|_{Y^{-s}}  \sum_{k=1}^5\|u_k\|_{X^s}\prod_{k=1;k\not=j}^5\|u_k\|_{X^1},
  \end{equation}
  We dyadically decompose each $\widetilde{u_k}$, and by symmetry in
  $u_1,\ldots,u_5$ it suffices to consider
  \begin{equation*}
    \Sigma:=\sum_{N_0\geq 0;\, N_1\geq\ldots\geq N_5\geq 0}
    \int_{\tau}\int_{ M} \prod_{j=0}^5 P_{N_j}\widetilde{u_j} \ dx dt
  \end{equation*}
  We split the sum $\Sigma=\Sigma_1+\Sigma_2$, where $\Sigma_1$ is
  defined by the constraint $\max\{N_0,N_2\}\sim N_1$. $\Sigma_1$ is
  the major contribution
  which can be handled by means of the Cauchy-Schwarz inequality and
  Corollary \ref{cor:tri-str-v2} precisely as in the proof of
  \cite[Proposition 4.1]{HTT10a}. The result is
  \begin{equation*}
    \Sigma_1\ls \| u_0\|_{Y^{-s}}  \|u_1\|_{X^s}\prod_{j=2}^5\|u_j\|_{X^1}.
  \end{equation*}
  In fact, in specific cases such as
    $M=\S^3$ there will be no further contribution because the product
    of five spherical harmonics of maximal degree $k$ can be developed
    into a series of spherical harmonics of maximal degree $5k$.
  In general, however, it remains to consider a minor contribution of lower
  order, which comes from the range where $\max\{N_0,N_2\}\ll N_1$ or
  $N_1\ll N_0$, and which we split $\Sigma_2=\Sigma_{21}+\Sigma_{22}$
  accordingly.  We have
  \begin{equation*}
    \Sigma_{21}\leq \sum_{ N_2\geq \ldots\geq
      N_5; \; N_1\gg N_0,N_2} \sum_{L\geq 0} \Big|\int_{\tau}
    I(N_0,\ldots,N_5,L)(t) dt \Big|,
  \end{equation*}
  where
  \begin{equation*}
    I(N_0,\ldots,N_5,L)(t)=\int_{ M} P_L\Big(\prod_{j=0}^2 P_{N_j}\widetilde{u_j}\Big)\prod_{j=3}^5 P_{N_j}\widetilde{u_j}  dx.
  \end{equation*}
  If $L\gs N_1$ we apply Lemma \ref{lem:decay} (recall that
  $N_3,N_4,N_5\ll L$), to deduce
  \begin{align*}
    &|I(N_0,\ldots,N_5,L)(t)|\\
    \ls{}& L^{-5} \|P_L(P_{N_0}u_0 P_{N_1}u_1 P_{N_2}u_2
    )(t)\|_{L^2(M)}\prod_{j=3}^5\|P_{N_j}u_j(t)\|_{L^2(M)},
  \end{align*}
  and H\"older's inequality and Lemma \ref{lem:crude} imply
  \begin{align*}
    \int_\tau |I(N_0,\ldots,N_5,L)(t)|dt \ls L^{-5}
    N_0^{\frac32}N_2^{\frac32}\prod_{j=0}^5\|P_{N_j}u_j\|_{L^\infty(\tau_0;L^2(M))},
  \end{align*}
  which implies
  \begin{equation}\label{eq:s21a}
    \sum_{L\gs N_1} \Big|\int_{\tau}
    I(N_0,\ldots,N_5,L)(t) dt \Big|
    \ls N_1^{-2}\prod_{j=0}^5\|P_{N_j}u_j\|_{V^2_\Delta}.
  \end{equation}
  On the other hand, if $L\ll N_1$ we apply Lemma \ref{lem:decay} (in
  this case $L,N_0,N_2\ll N_1$), to deduce
  \begin{align*}
    &|I(N_0,\ldots,N_5,L)(t)|\\
    \ls{}& N_1^{-5} \prod_{j=0}^2\|P_{N_j}u_j(t)\|_{L^2(M)}
    \|P_L(P_{N_3}u_3 P_{N_4}u_4 P_{N_5}u_5)(t)\|_{L^2(M)}
  \end{align*}
  and H\"older's inequality and Lemma \ref{lem:crude} imply
  \begin{align*}
    \int_\tau |I(N_0,\ldots,N_5,L)(t)|dt \ls N_1^{-5}
    N_4^{\frac32}N_5^{\frac32}\prod_{j=0}^5\|P_{N_j}u_j\|_{L^\infty(\tau_0;L^2(M))}.
  \end{align*}
  which together with \eqref{eq:s21a} gives
  \begin{align*}
    \sum_{L\geq 0} \Big|\int_{\tau} I(N_0,\ldots,N_5,L)(t) dt \Big|
    \ls N_1^{-1}\prod_{j=0}^5\|P_{N_j}u_j\|_{V^2_\Delta}.
  \end{align*}
  Dyadic summation easily yields
  \begin{equation*}
    \Sigma_{21}\ls  \|u_0\|_{Y^{-s}}\|u_1\|_{X^{s}}\prod_{j=1}^5\|u_j\|_{X^1}.
  \end{equation*}
  The contribution of $\Sigma_{22}$, where $N_0\gg N_1\geq N_2\geq
  \ldots\geq N_5$ can be treated in the same way by switching the
  roles of $N_1$ and $N_0$.
\end{proof}

The proof of Theorem \ref{thm:main-tech} -- based on Proposition
\ref{prop:six-est} and the contraction mapping principle -- is
standard and can be concluded as in \cite[Section 4]{HTT10a}, cp. also
the references therein.
\appendix

\section{Proof of Lemma \ref{lem:lp}}\label{asect:p}
For the sake of completeness, we include a proof of Lemma \ref{lem:lp}
here. This result has been proved in the case $J=[1,N]$ and
$\mu_n^2=n^2$ in \cite[Section 4]{B89} and stated for general $J$ in
the case $p=6$ in \cite{BGT07}.  More precisely, we describe here the
necessary modifications with respect to Bourgain's original work
\cite[Section 4]{B89}. We closely follow the presentation in
\cite[Section 4]{B89}. For
  this reason we work in the $1$-periodic (instead of the
  $2\pi$-periodic) setup here. However, note that we replaced $\delta
  N^{\frac12}$ with $\lambda$. We also
refer the reader to \cite[Section 3]{B93a}, and to the book \cite{V81}
for more details on the circle method of Hardy and Littlewood.

First of all, we switch from $t$ to $-t$, and reduce the estimate
\eqref{lem:lp} for general $\alpha\in \N_0$ to the case $\alpha=0$.
The latter simply follows by dilating time by the factor $16$ and
applying the result for $\alpha=0$ on $[0,2\pi]$ to a modified
sequence and the translated and dilated interval $4J+\alpha$. From now
on we will assume that $\mu_n^2=n^2$.

Let $\N$ denote the set of positive integers, and let $J=[b,b+N]$,
$b,N\in \N$. As in \cite[Section 3]{B93a} we choose a sequence
$\sigma$ satisfying
\begin{enumerate}
\item For all $n\in \Z$: $0\leq \sigma_n\leq 1$; for all $n \in J$:
  $\sigma_n=1$ for $n \in J$; for all $n$ such that $n< b-N$ or
  $n>b+2N$: $\sigma_n=0$.
\item The sequence $(\sigma_{n+1}-\sigma_n)$ is bounded by $N^{-1}$
  and has variation bounded by $N^{-1}$.
\end{enumerate}
Let $p>4$, and $0<\eps\ll 1$ such that $p-\eps>4$. Our aim is to prove
the distributional inequality
\begin{equation}\label{eq:distr}
  \sup_{b \in \N}\Big|\Big\{t \in [0,1]: \Big|\sum_{n\in \Z}
  c_n\sigma_n e^{2\pi it n^2}\Big|>\lambda\Big \}\Big|\leq C_\eps N^{\eps/2}\lambda^{-4-\eps},
\end{equation}
for all $c_n$ such that $\sum_{n}|c_n|^2=1$, which implies
\eqref{eq:lp}, because the set is empty if $\lambda\gg N^{1/2}$.

{\it a)} There exists $c_\eps>0$ such that
\[
\sup_{k\in \N, b \in \N_0}\#\{(n_1,n_2)\in \N^2: n_1,n_2\leq N;\quad
n_1(n_2+b)=k\}\leq c_\eps N^{\frac{\eps}{400}}.
\]
In the case $0\leq b\leq 10N^2$ this follows from the standard bound
on the number of divisors function, see \cite[Theorem 315]{HW79}.
Otherwise, the set contains at most one element.  With this ingredient
one can easily modify the argument in \cite[formulas (1.3)-(1.6)]{B89}
to deduce
\begin{equation}\label{eq:l4}
  \Big\|\sum_{n\in \Z} c_n\sigma_n e^{2\pi it n^2}\Big\|_{L^4_t(0,1)}\ls {}
  N^{\frac{\eps}{400}}\Big(\sum_{n\in \Z} |c_n|^2\Big)^{\frac12}.
\end{equation}
This bound implies \eqref{eq:distr} for $\lambda\ls {} N^{1/2-\nu/4}$
for $\nu=1/100$, so it remains to prove \eqref{eq:distr} in the case
$N^{1/2-\nu/4} \ll \lambda\leq N^{1/2}$.

{\it b)} It holds
\begin{equation}\label{eq:lem318-alt}
  \Big|\sum_{n \in \Z} \sigma_n e^{2\pi it n^2}\Big|\ls {} q^{-1/2}(|t-a/q|+N^{-2})^{-1/2},
\end{equation}
for any $1\leq a < q< N$, $\gcd(a,q)=1$ and $|t-a/q|<(qN)^{-1}$.  The
claim \eqref{eq:lem318-alt} follows from
\[
\sum_{n\in \Z} \sigma_n e^{2\pi it n^2}=e^{2\pi it b^2}\sum_{m \in \Z
} \sigma_{b+m} e^{4\pi i bt }e^{2\pi it m^2}
\]
and \cite[Lemma 3.18]{B93a} with $x=2bt$.  Estimate
\eqref{eq:lem318-alt} replaces \cite[formula (4.10)]{B89}, with
\[f(t)=\sum_{n \in \Z} \sigma_n e^{2\pi it n^2}.\]

{\it c)} We define the major arcs $\mathcal{M}$ to be the disjoint
union of the sets
\[\mathcal{M}(q,a)=\{t \in [0,1]: |t-a/q|\leq N^{\nu-2}\},\]
for any $1\leq a\leq q\leq N^\nu$, $\gcd(a,q)=1$.  Let $t \in
[0,1]\setminus \mathcal{M}$. By Dirichlet's Lemma there exists a
reduced fraction $a/q$ with $1\leq a\leq q\leq N^{2-\nu}$ such that
$|t-a/q|\leq N^{\nu-2}$, and since $t \notin\mathcal{M}$ it must be
$q>N^\nu$ and \eqref{eq:lem318-alt} implies
\begin{equation}\label{eq:outsidemj}
  |f(t)|\leq N^{1-\nu/2},
\end{equation} which replaces \cite[formula (4.6)]{B89}.

{\it d)} In order to prove \eqref{eq:distr}, it therefore suffices to
prove a bound on the number $R$ of $N^{-2}$-separated points
$t_1,\ldots,t_R \in [0,1]$ where
\[
\Big|\sum_{n \in \Z} c_n\sigma_n e^{2\pi i n^2t_r} \Big|>\lambda.
\]
We recall that $N^{1/2-\nu/4} \ll \lambda\leq N^{1/2}$. As in
\cite{B89} we obtain
\[\sum_{1\leq r,r'\leq R}\Big|\sum_{n \in \Z} \sigma_n e^{2\pi i
  n^2(t_r-t_{r'})}\Big|>\lambda^2 R^2.\] For fixed $\gamma>2$ this
estimate and H\"older's inequality yield
\[
\sum_{1\leq r,r'\leq R} |f(t_r-t_{r'})|^{\gamma}>\lambda^{2\gamma}
R^2,
\]
which replaces \cite[formula (4.13) with $\lambda=\delta
N^{\frac12}$]{B89}. From here, the arguments in \cite[pp.
305--307]{B89} apply verbatim and show that
\[R\ls {} \lambda^{-4-\eps}N^{2+\frac{\eps}{2}}.\] Because of the
$N^{-2}$-separation property of the points $t_1,\ldots,t_R$ this
implies
\begin{equation*}
  \Big|\Big\{t \in [0,1]: \Big|\sum_{n\in \Z} \sigma_n c_n e^{2\pi it n^2}\Big|>\lambda\Big\}\Big|\ls {} N^{-2} \lambda^{-4-\eps}N^{2+\frac{\eps}{2}},
\end{equation*}
which gives \eqref{eq:distr}.

\subsection*{Acknowledgments}
I am indebted to Daniel Tataru and Nikolay Tzvetkov for introducing me
to this circle of problems and for helpful and stimulating
discussions. I would also like to thank Christoph Thiele for a helpful
conversation about \cite{B89}.

\bibliographystyle{amsplain} \bibliography{nls-refs}\label{sect:refs}

\end{document}